\documentclass[12pt,thmsb]{article}
\usepackage{graphicx,epsfig}
\usepackage{amsmath,amsbsy}
\usepackage{amsfonts}
\usepackage{amssymb}
\usepackage{amsmath}
\usepackage[onehalfspacing]{setspace}
\usepackage{graphicx}

\newtheorem{theorem}{Theorem}
\newtheorem{corollary}{Corollary}
\newtheorem{lemma}{Lemma}
\newtheorem{definition}{Definition}
\newtheorem{remark}{Remark}

\newtheorem{proposition}{Proposition}

\def\1{\hbox{\rm 1\hskip -3pt I}}

\input{tcilatex}

\begin{document}

\author{Fr\'{e}d\'{e}ric \textsc{Ferraty} \thanks{%
Universit\'{e} Paul Sabatier, Toulouse 3 and Universit\'{e} Toulouse le
Mirail, Toulouse 2} , Andr\'{e} \textsc{Mas} \thanks{
Universit\'{e} Montpellier 2} \thanks{%
Corresponding author: Laboratoire de Probabilit\'{e}s et Statistique CC051,
Universit\'{e} Montpellier 2, Place Eug\`{e}ne Bataillon, 34095 Montpellier,
France, mas@math.univ-montp2.fr}, and Philippe \textsc{Vieu}\thanks{
Universit\'{e} Paul Sabatier, Toulouse 3}}
\title{Nonparametric regression on functional data: inference and practical
aspects}
\maketitle

\begin{abstract}
We consider the problem of predicting a real random variable from a
functional explanatory variable. The problem is attacked by mean of
nonparametric kernel approach which has been recently adapted to this
functional context. We derive theoretical results by giving a deep
asymptotic study of the behaviour of the estimate, including mean squared
convergence (with rates and precise evaluation of the constant terms) as
well as asymptotic distribution. Practical use of these results are relying
on the ability to estimate these constants. Some perspectives in this
direction are discussed. In particular a functional version of wild
bootstrapping ideas is proposed and used both on simulated and real
functional datasets.
\end{abstract}



\newpage 

\section{Introduction}

\label{Intro} 

Functional data are more and more frequently involved in statistical
problems. Developping statistical methods in this special framework has been
popularized during the last few years, particularly with the monograph by
Ramsay \& Silverman (2005). More recently, new developments have been
carried out in order to propose nonparametric statistical methods for
dealing with such functional data (see Ferraty \& Vieu, 2006, for large
discussion and references). These methods are also called doubly infinite
dimensional (see Ferraty \& Vieu, 2003). Indeed these methods deal with
infinite-dimensional (i.e. functional) data and with a statistical model
which depends on an infinite-dimensional unknown object (i.e. a
nonparametric model). This double infinite framework motivates the
appellation of Nonparametric Functional Statistics for such kind of methods.
Our paper is centered on the functional regression model :%
\begin{equation}
Y=r(\mathcal{X})+error,
\end{equation}
where $Y$ is a real random variable, $\mathcal{X}$ is a functional random
variable (that is, $\mathcal{X}$ takes values in some possibly
infinite-dimensional space) and where the statistical model assumes only
smoothness restriction on the functional operator $r$. At this point, it
worth noting that the operator $r$ is not constrained to be linear. This is
a Functional Nonparametric Regression model (see Section \ref{Notations} for
deeper presentation).

The aim of this paper is to extend in several directions the current
knowledges about functional nonparametric regression estimates presented in
Section \ref{Notations}. In Section \ref{MSE} we give asymptotic mean
squared expansions, while in Section \ref{AsNorm} the limiting distribution
is derived. The main novelty/difficuly along the statement of these results
relies on the exact calculation of the leading terms in the asymptotic
expressions. Section \ref{exSBP} points out how such results can be used
when the functional variable belongs to standard families of continuous time
process. The accuracy of our asymptotic results leads to interesting
perspectives from a practical point of view: minimizing mean squared errors
can govern automatic bandwidth selection procedure while the limiting
distribution of the error is a useful tool for building confidence bands. To
this end, we propose in Section \ref{ComputFeatures} a functional version of
the wild bootstrap procedure, and we use it, both on simulated and on real
functional datasets, to get some automatic rule for choosing the bandwidth.
The concluding section \ref{Conc} contains some important open questions
which emerge naturally from the theoretical results given in this paper,
such as the theoretical study of the accuracy of the functional wild
bootstrap procedure used in our applications.



\section{Kernel nonparametric functional regression}

\label{Notations} 

\subsection{The model}

The model is defined in the following way. Assume that $\left( Y_{i},%
\mathcal{X}_{i}\right) $ is a sample of $n$ i.i.d. pairs of random
variables. The random variables $Y_{i}$ are real and the $\mathcal{X}_{i}$'s
are random elements with values in a functional space $\mathcal{E}$. In all
the sequel we will take for $\mathcal{E}$ a separable Banach space endowed
with a norm $\left\Vert \cdot\right\Vert $. This setting is quite general
since it contains the space of continuous functions, $L^{p}$ spaces as well
as more complicated spaces like Sobolev or Besov spaces. Separability avoids
measurability problems for the random variables $\mathcal{X}_{i}$'s. The
model is classically written : 
\begin{equation*}
Y_{i}=r\left( \mathcal{X}_{i}\right) +\varepsilon_{i}\quad i=1,...,n,
\end{equation*}
where $r$ is the regression function mapping $\mathcal{E}$ onto $\mathbb{R}$
and the $\varepsilon_{i}$'s are such that for all $i$, $E(\varepsilon _{i}|%
\mathcal{X}_{i})=0$ and $E(\varepsilon_{i}^{2}|\mathcal{X}%
_{i})=\sigma_{\varepsilon}^{2}(\mathcal{X})<\infty$.


\subsection{The estimate}


Estimating $r$ is a crucial issue in particular for predicting the value of
the response given a new explanatory functional variable $\mathcal{X}_{n+1}$%
. However, it is also a very delicate task because $r$ is a nonlinear
operator (from $\mathcal{E}$ into $\mathbb{R}$) for which functional linear
statistical methods were not planned. To provide a consistent procedure to
estimate the nonlinear regression operator $r$, we propose to adapt the
classical finite dimensional Nadaraya-Watson estimate to our functional
model. We set 
\begin{equation*}
\widehat{r}\left( \chi\right) =\frac{\sum_{k=1}^{n}Y_{k}K\left(
h^{-1}\left\Vert \mathcal{X}_{k}-\chi\right\Vert \right) }{%
\sum_{k=1}^{n}K\left( h^{-1}\left\Vert \mathcal{X}_{k}-\chi\right\Vert
\right) }.
\end{equation*}

\noindent Several asymptotic properties of this estimate were obtained
recently. It turns out that the existing literature adresses either the
statement of upper bounds of the rates of convergence without specification
of the exact constants (see Chapter 6 in Ferraty \& Vieu, 2006), or abstract
expressions of these constants which are unusable in practice (as for
instance in the recent work by Masry, 2005, which has been published during
the reviewing process of this paper). Our aim in this paper is to give bias,
variance, means square errors and asymptotic distribution of the functional
kernel regression estimate with exact computation of all the constants (see
Section \ref{Theorie}). We will focus on practical purposes in Section \ref%
{ComputFeatures}.

Several assumptions will be made later on the kernel $K$ and on the
bandwidth $h$. Remind that in a finite-dimensional setting pointwise mean
squared error (at $\chi$) of the estimate depends on the evaluation of the
density (at $\chi$) w.r.t. Lebesgue's measure and on the derivatives of this
density. We refer to Schuster (1972) for an historical result about this
topic. On infinite-dimensional spaces, there is no measure universally
accepted (as the Lebesgue one in the finite-dimensional case) and there is
need for developping a \textquotedblleft free-density\textquotedblright\
approach. As discussed along Section \ref{exSBP} the problem of introducing
a density for $\mathcal{X}$ is shifted to considerations on the measure of
small balls with respect to the probability of $X$.

\subsection{Assumptions and notations}

\label{Assump} Only pointwise convergence will be considered in the
forthcoming theoretical results. In all the following, $\chi$ is a fixed
element of the functional space $\mathcal{E}$. Let $\varphi$ be the real
valued function defined as 
\begin{equation*}
\varphi\left( s\right) =E\left[ \left( r\left( \mathcal{X}\right) -r\left(
\chi\right) \right) |\left\Vert \mathcal{X}-\chi\right\Vert =s\right] ,
\end{equation*}
and $F$ be the c.d.f. of the random variable $\left\Vert \mathcal{X}%
-\chi\right\Vert $: 
\begin{equation*}
F\left( t\right) =P\left( \left\Vert \mathcal{X}-\chi\right\Vert \leq
t\right) .
\end{equation*}
Note that the crucial functions $\varphi$ and $F$ depends implicitely on $%
\chi.$ Consequently we should rather note them by $\varphi_{\chi}$ and $%
F_{\chi}$ but, as $\chi$ is fixed, we drop this index once and for all.
Similarly, we will use in the remaining the notation $\sigma_{%
\varepsilon}^{2}$ instead of $\sigma_{\varepsilon}^{2}(\mathcal{X})$. Let us
consider now the following assumptions.

\textbf{H0 :} $r$ and $\sigma^{2}_{\varepsilon}$ are continuous in a
neighborhood of $\chi$, and $F(0)=0$.

\textbf{H1} \textbf{:} $\varphi^{\prime}(0)$ exists.

\textbf{H2 :} The bandwidth $h$ satisfies $\lim_{n\rightarrow\infty}h=0$ and 
$\lim_{n\rightarrow\infty}nF(h)=\infty$, while the kernel $K$ is supported
on $[0,1]$, has a continuous derivative on $[0,1)$, $K^{\prime}(s)\leq0$ and 
$K\left( 1\right) >0$.

Assumptions \textbf{H0} and \textbf{H2} are clearly unrestrictive, since
they are the same as those classically used in the finite-dimensional
setting. Much more should be said on assumption \textbf{H1}. Note first
that, obviously, $\varphi\left( 0\right) =0$. It is worth noting that,
whereas we could expect assumptions on the local regularity of $r$ (as in
the finite-dimensional case), hypothesis \textbf{H1} skips over that point
and avoids to go into formal considerations on differential calculus on
Banach spaces. To fix the ideas, if we assumed differentiability of $r$, we
would get by Taylor's expansion that 
\begin{equation*}
r\left( \mathcal{X}\right) -r\left( \chi\right) =\left\langle r^{\prime
}\left( \chi\right) ,\mathcal{X}-\chi\right\rangle +o\left( \left\Vert 
\mathcal{X}-\chi\right\Vert \right) ,
\end{equation*}
where $r^{\prime}\left( \chi\right) \in\mathcal{E}^{\ast}$, $\mathcal{E}%
^{\ast}$ being the conjugate space of $\mathcal{E}$ and $\left\langle
\cdot,\cdot\right\rangle $ being the duality bracket between $\mathcal{E}$
and $\mathcal{E}^{\ast}.$ In this context, a non trivial link would appear
between $\varphi^{\prime}(0)$ and $r^{\prime}(\chi)$ through the following
relation: 
\begin{equation*}
\lim_{s\rightarrow0}E\left[ \left\langle r^{\prime}\left( \chi\right) ,\frac{%
\mathcal{X}-\chi}{\left\Vert \mathcal{X}-\chi\right\Vert }\right\rangle
|\left\Vert \mathcal{X}-\chi\right\Vert =s\right] =\varphi^{\prime}\left(
0\right) .
\end{equation*}
Indeed, even if the link between the existency of $r^{\prime}(\chi)$ and of $%
\varphi^{\prime}(0)$ is strong, one can build counter-examples for their non
equivalence (these counter-examples are available on request but they are
out of the main scope of this paper). In the perspective of estimating the
constants given in Theorem \ref{TheoremMSE}, it will be easier to estimate $%
\varphi^{\prime}(0)$ (for instance by using $\widehat{r}$) than the operator 
$r^{\prime}(\chi)$. Therefore, we prefer to express computations by mean of $%
\varphi^{\prime}(0)$ instead of $r^{\prime}(\chi)$. This has the additional
advantage to produce more readable writings. Consequently, the
differentiability of $r$ is not needed.

Let us now introduce the function $\tau_{h}$ defined for all $s\in\left[ 0,1%
\right] $ as: 
\begin{equation*}
\tau_{h}\left( s\right) =\frac{F\left( hs\right) }{F\left( h\right) }%
=P\left( \left\Vert \mathcal{X}-\chi\right\Vert \leq hs|\left\Vert \mathcal{X%
}-\chi\right\Vert \leq h\right) ,
\end{equation*}
for which the following assumption is made:

\textbf{H3 :} For all $s\in\left[ 0,1\right] ,$ $\tau_{h}\left( s\right)
\rightarrow\tau_{0}\left( s\right) $ as $h \rightarrow0$.

\noindent Note that the function $\tau_{h}$ is increasing for all $h.$ The
measurable (as the pointwise limit of the sequence of measurable functions $%
\tau_{h}$) mapping $\mathbf{\tau}_{0}$ is non decreasing. Let us finally
mention that this function $\mathbf{\tau}_{0}$ will play a key role in our
methodology, in particular when we will have to compute the exact constant
terms involved in our asymptotic expansions. For the sake of clarity, the
following proposition (whose a short proof will be given in the Appendix)
will explicit the function $\tau_{0}$ for various cases. By $1_{\left] 0,1%
\right] }(\cdot)$ we denote the indicator function on the set $\left] 0,1%
\right] $ and $\delta_{1}\left( \cdot\right) $ stands for the Dirac mass at $%
1$.

\begin{proposition}
\label{prop1}~\newline
i) If $F\left( s\right) \sim Cs^{\gamma}$ for some $\gamma>0$ then $%
\tau_{0}(s)=s^{\gamma}$.\newline
ii) If $F\left( s\right) \sim Cs^{\gamma}\left\vert \ln s\right\vert
^{\kappa}$ with $\gamma>0$ and $\kappa>0$ then $\tau_{0}\left( s\right)
=s^{\gamma}.$ \newline
iii) If $F\left( s\right) \sim C_{1}s^{\gamma}\exp\left( -C_{2}/s^{p}\right) 
$ for some $p>0$ and some $\gamma>0$ then $\tau_{0}\left( s\right) =\delta
_{1}\left( s\right) .$ \newline
iv) If $F\left( s\right) \sim C/\left\vert \ln s\right\vert $ then $%
\tau_{0}\left( s\right) =1_{\left] 0,1\right] }(s)$.
\end{proposition}

A deeper discussion linking the above behavior of $F$ with small ball
probabilities notions will be given in Section \ref{exSBP}.


\section{Asymptotic study}

\label{Theorie} 

In both following subsections we will state some asymptotic properties
(respectively mean squared asymptotic evaluation and asymptotic normality)
for the functional kernel regression estimate $\hat{r}$.

It is worth noting that all the results below can be seen as extensions to
functional data of several ones already existing in the finite-dimensional
case (the literature is quite extensive in this field and the reader will
find in Sarda \& Vieu (2000) deep results as well as a large scope of
references). With other words, our technique for proving both Theorem \ref%
{TheoremMSE} and Theorem \ref{TheoremAsNorm} is also adapted to the scalar
or vector regression model since the abstract space $\mathcal{E}$ can be of
finite dimension (even, of course, if our main goal is to treat
infinite-dimensional cases). Moreover, it turns out that the transposition
to finite-dimensional situations of our key conditions (see discussion in
Section \ref{exSBP} below) becomes (in some sense) less restrictive than
what is usually assumed. With other words, the result of Theorem \ref%
{TheoremMSE} and Theorem \ref{TheoremAsNorm} can be directly applied to
finite-dimensional settings, and will extend the results existing in this
field (see again Sarda \& Vieu, 2000) to situation when the density of the
corresponding scalar or multivariate variable does not exist or has all its
successive derivatives vanishing at point $\chi$ (see discussion in Section %
\ref{DimFinie}).

All along this section we assume that assumptions \textbf{H0-H3} hold. Let
us first introduce the following notations:

\begin{align*}
M_{0} & =\left( K\left( 1\right) -\int_{0}^{1}\left( sK\left( s\right)
\right) ^{\prime}\tau_{0}\left( s\right) ds\right) , \\
M_{1} & =\left( K\left( 1\right) -\int_{0}^{1}K^{\prime}\left( s\right)
\tau_{0}\left( s\right) ds\right) , \\
M_{2} & =\left( K^{2}\left( 1\right) -\int_{0}^{1}\left( K^{2}\right)
^{\prime}\left( s\right) \tau_{0}\left( s\right) ds\right) .
\end{align*}


\subsection{Mean Squared Convergence}

\label{MSE} 
The following result gives asymptotic evaluation of the mean squared errors
of our estimate. The asymptotic mean squared errors have a standard convex
shape, with large bias when the bandwidth $h$ increases and large variance
when $h$ decays to zero. We refer to the Appendix for the proof of Theorem %
\ref{TheoremMSE}.

\begin{theorem}
\label{TheoremMSE} When \textbf{H0-H3} hold, we have the following
asymptotic developments :%
\begin{equation}
E\widehat{r}\left( \chi\right) -r\left( \chi\right) =\varphi^{\prime }\left(
0\right) \frac{M_{0}}{M_{1}}h+O\left( \left( nF\left( h\right) \right)
^{-1}\right) +o\left( h\right) ,  \label{B}
\end{equation}
and%
\begin{equation}
Var\left( \widehat{r}\left( \chi\right) \right) =\frac{1}{nF\left( h\right) }%
\frac{M_{2}}{M_{1}^{2}}\sigma_{\varepsilon}^{2}+o\left( \frac {1}{nF\left(
h\right) }\right) .  \label{V}
\end{equation}
\end{theorem}


\subsection{Asymptotic Normality}

\label{AsNorm} 

Let us denote the leading bias term by: 
\begin{equation*}
B_{n}=\varphi^{\prime}\left( 0\right) \frac{M_{0}}{M_{1}}h.
\end{equation*}
Before giving the asymptotic normality, one has to be sure that the leading
bias term does not vanish. This is the reason why we introduce the following
additional assumption:

\textbf{\ H4 :} $\varphi^{\prime}(0)\neq0$ and $M_{0}>0$.\newline
The first part of assumption \textbf{H4} is very close to what is assumed in
standard finite-dimensional literature. It forces the nonlinear operator $r$
not to be too smooth (for instance, if $r$ is Lipschitz of order $\beta>1$,
then $\varphi^{\prime}(0)=0$). The second part of assumption \textbf{H4} is
specific to the infinite-dimensional setting, and the next Proposition \ref%
{Pr2} will show that this condition is general enough to be satisfied in
some standard situations. This proposition will be proved in the appendix.

\begin{proposition}
\label{Pr2}~\newline
i) If $\tau_{0}(s)\neq1_{\left] 0,1\right] }(s)$ and $\tau_{0}$ is
continuously differentiable on $(0,1)$, then $M_{0}>0$ for any kernel $K$
satisfying \textbf{H2}.\newline
ii) If $\tau_{0}(s)=\delta_{1}(s)$, then $M_{0}>0$ for any kernel $K$
satisfying \textbf{H2}.
\end{proposition}

To emphasize the interest of these results, they should be combined with
those of Proposition \ref{prop1}. Note that the result $i)$ includes the
well-known family of processes for which $\tau_{0}(s)=s^{\gamma}$ (see
Proposition \ref{prop1}-$i$), that is those whose distributions admit
fractal dimensions (see Section \ref{fractal}). The second case when $%
\tau_{0}(s)=\delta_{1}(s)$ (for which a particular case is given in
Proposition \ref{prop1}-$iii$) corresponds to nonsmooth processes (see
Section \ref{nonsmooth}). These two cases cover a large number of
situations. However, if a more general function $\tau_{0}$ has to be used,
one can make additional hypotheses on the kernel $K$. In particular, if $K$
is such that $\forall s\in\lbrack0,1]$, $(sK(s))^{\prime}>0$ and $%
\tau_{0}(s)\neq1_{\left] 0,1\right] }(s)$, then $M_{0}>0$, which covers the
case of the uniform kernel. More complicated kernel functions $K$ would lead
to more technical assumptions linking $K$ with $\tau_{0}$. It is out of
purpose to give these tedious details (available on request) but let us just
note that the key restriction is the condition $\tau_{0}(s)\neq1_{\left] 0,1%
\right] }(s)$ (else we have $M_{0}=0$).

Moreover, since the rate of convergence depends on the function $F\left(
h\right) $ and for producing a reasonably usable asymptotic distribution it
is worth having some estimate of this function. The most natural is its
empirical counterpart: 
\begin{equation*}
\widehat{F}\left( h\right) =\frac{\#\left( i:\left\Vert \mathcal{X}%
_{i}-\chi\right\Vert \leq h\right) }{n}.
\end{equation*}

\noindent The pointwise asymptotic gaussian distribution for the functional
nonparametric regression estimate is given in Theorem \ref{TheoremAsNorm}
below which will be proved in the appendix. Note that the symbol $%
\hookrightarrow$ stands for ''convergence in distribution''.

\begin{theorem}
\label{TheoremAsNorm} When \textbf{H0-H4} hold, we have 
\begin{equation*}
\sqrt{n\widehat{F}\left( h\right) }\left( \widehat{r}\left( \chi\right)
-r\left( \chi\right) -B_{n}\right) \frac{M_{1}}{\sqrt{M_{2}\sigma
_{\varepsilon}^{2}}}\hookrightarrow\mathcal{N}\left( 0,1\right) .
\end{equation*}
\end{theorem}

\noindent A simpler version of this result is stated in Corollary \ref{coro1}
below whose proof is obvious. The key-idea relies in introducing the
following additional assumption:

\textbf{H5 :} $\lim_{n\rightarrow\infty}h\sqrt{n\, F(h)}\ =\ 0$

\noindent which allows to cancel the bias term.

\begin{corollary}
\label{coro1} When \textbf{H0-H5} hold, we have 
\begin{equation*}
\sqrt{n\widehat{F}\left( h\right) }\left( \widehat{r}\left( \chi\right)
-r\left( \chi\right) \right) \frac{M_{1}}{\sigma_{\varepsilon}\sqrt{M_{2}}}%
\hookrightarrow\mathcal{N}\left( 0,1\right) .
\end{equation*}
\end{corollary}

\noindent In practice, the constants involved in Corollary \ref{coro1} need
to be estimated. In order to compute explicitely both constants $M_{1}$ and $%
M_{2}$, one may consider the simple uniform kernel and get easily the
following result:

\begin{corollary}
Under assumptions of Corrollary \ref{coro1}, if $K(.)=1_{[0,1]}(.)$ and if $%
\widehat{\sigma}_{\varepsilon}^{2}$ is a consistent estimator of $%
\sigma_{\varepsilon}^{2}$, then we have: 
\begin{equation*}
\sqrt{\frac{n\widehat{F}\left( h\right) }{\widehat{\sigma}_{\varepsilon}^{2}}%
}\left( \widehat{r}\left( \chi\right) -r\left( \chi\right) \right)
\hookrightarrow\mathcal{N}\left( 0,1\right) .
\end{equation*}
\end{corollary}

\noindent There are many possibilities for constructing a consistent
conditional variance estimate. One among all the possibilities consists in
writing that 
\begin{align*}
\sigma^{2}_{\varepsilon}(\chi) & = E((Y-r(\mathcal{X}))^{2}|\mathcal{X}%
=\chi), \\
& = E(Y^{2}|\mathcal{X}=\chi) - (E(Y|\mathcal{X}=\chi))^{2},
\end{align*}
and, by estimating each conditional expectation with the functional kernel
regression technique.


\section{Some Examples of small ball probabilities}

\label{exSBP}

\label{Examples} 

The distribution function $F$ plays a prominent role in our methodology.
This appears clearly in our conditions (through the function $\tau_{0}$) and
in the rates of convergence of our estimate (through the asymptotic behavior
of the quantity $n\,F(h)$). More precisely, the behaviour of $F$ around $0$
turns out to be of first importance. In other words, the small ball
probabilities of the underlying functional variable $\mathcal{X}$ will be
determining. In order to illustrate our ideas and to connect with existing
probabilistic knowledges in this field, let us now just discuss how $F$ (and
hence $\tau_{0}$) behave for different usual examples of processes $\mathcal{%
X}$ valued in an infinite-dimensional space.

\subsection{Nonsmooth processes}

\label{nonsmooth}

Calculation of the quantity $P\left( \left\Vert \mathcal{X}-\chi\right\Vert
<s\right) $ for \textquotedblleft small\textquotedblright\ $s$ (i.e. for $s$
tending to zero) and for a fixed $\chi$ is known as a "small ball problem"
in probability theory. This problem is unfortunately solved for very few
random variables (or processes) $\mathcal{X},$ even when $\chi=0.$ In
certain functional spaces, taking $\chi\neq0$ yield considerable
difficulties that may not be overcome. Authors usually focus on gaussian
random elements. We refer to Li \& Shao (2001) for a survey on the main
results on small ball probability. If $\mathcal{X}$ is a gaussian random
element on the separable Banach space $\mathcal{E}$ and if $\chi$ belongs to
the reproducing kernel Hilbert space associated with $\mathcal{X}$, then the
following well-known result holds: 
\begin{equation}
P\left( \left\Vert \mathcal{X}-\chi\right\Vert <s\right) \sim C_{\chi
}P\left( \left\Vert \mathcal{X}\right\Vert <s\right) ,\quad\mathrm{as}\text{ 
}s\rightarrow0.  \label{pb}
\end{equation}
So, the small ball problem at any point $\chi$ may be shifted to a small
ball problem at $0.$ Moreover, (\ref{pb}) can be precised in a few
situations. For instance, Mayer-Wolf \& Zeitouni (1993) investigate the case
when $\mathcal{X}$ is a one-dimensional diffusion process and $\chi$
satisfies some conditions (see Mayer-Wolf \& Zeitouni, 1993, p15). They also
briefly mention the non gaussian case (see Mayer-Wolf \& Zeitouni, 1993,
Remark 3, p19) but many other authors have considered different settings
(see Ferraty et \emph{al.}, 2005, for a large discussion and references
therein). As far as we know, the results which are available in the
literature are basically all of the form: 
\begin{equation}
P\left( \left\Vert \mathcal{X}-\chi\right\Vert <s\right) \sim c_{\chi
}\,s^{-\alpha}\exp\left( -C/s^{\beta}\right) ,  \label{Ex1}
\end{equation}
where $\alpha,\beta,c_{\chi}$ and $C$ are positive constants and $\left\Vert
\cdot\right\Vert $ may be a $\sup$, a $L^{p}$, a Besov norm ... The next
remark is a direct consequence of Proposition \ref{prop1}. It proves that
non-smooth processes may satisfy the assumptions needed to get the
asymptotic expansions of previous sections.

\begin{remark}
\label{Rem1} In the case of "non-smooth" processes defined by (\ref{Ex1}) we
have $\tau_{0}\left( s\right) =\delta_{1}(s)$. In addition, condition $%
n\,F\left( h\right) \rightarrow+\infty$ (in \textbf{H2}) is checked as soon
as $h=A/\left( \log n\right) ^{1/\beta}$ for $A$ large enough.
\end{remark}

\subsection{Fractal (or geometric) processes}

\label{fractal}

Another family of infinite dimensional processes is the class of fractal
processes for which the small ball probabilities are of the form 
\begin{equation}
P\left( \left\Vert \mathcal{X}-\chi\right\Vert <s\right) \sim c^{\prime
}_{\chi}\, s^{\gamma},  \label{Ex2}
\end{equation}
where $c^{\prime}_{\chi}$ and $\gamma$ are once again positive constants.
Like above, it is elementary to get the following result from Proposition %
\ref{prop1}.

\begin{remark}
\label{Rem2} Under (\ref{Ex2}), we have $\tau_{0}(s)=s^{\gamma}$ while the
condition $n\, F\left( h\right) \rightarrow+\infty$ (in \textbf{H2}) is
satisfied as soon as $h=A n^{-B}$ for $B$ small enough.
\end{remark}

\subsection{Back to the finite dimensional setting}

\label{DimFinie}

Finally, it is important to note that a special case of fractal processs is
given by the usual multivariate case (that is, by the case when $\mathcal{E}=%
\mathbb{R}^{p}$). The following result is obvious for the uniform norm on $%
\mathbb{R}^{p}$ and extends directly to any norm, since all of them are
equivalent in finite dimension.

\begin{remark}
\label{Rem3} If $\mathcal{E}=\mathbb{R}^{p}$, then any random variable $%
\mathcal{X}$ on $\mathbb{R}^{p}$ which has a finite and non zero density
function at point $\chi$ satisfies (\ref{Ex2}) with $\gamma=p$.
\end{remark}

From Remarks \ref{Rem2} and \ref{Rem3}, it is clear that all the results of
Section \ref{Theorie} apply in a finite dimensional setting. Besides, the
assumptions needed for Theorems \ref{TheoremMSE} and \ref{TheoremAsNorm} are
weaker than those described in Remark \ref{Rem3} since there is no need to
assume the existence of a density for $\mathcal{X}$. In this sense, our
results extend the standard multivariate literature (see the discussion at
the beginning of Section \ref{Theorie}).



\section{Perspectives on bandwidth choice}

\label{ComputFeatures} 


\subsection{Introduction}

\label{IntroBoot} 

The asymptotic results presented in the previous Section \ref{Theorie} are
particularly appealing because, in addition to the specification of the
rates of convergence, the exact constants involved in the leading terms of
each result are precised. This is particularly interesting in practice. Let
us focus now on the mean squared errors expansion given in Section \ref{MSE}%
. In fact, Theorem \ref{TheoremMSE} could give clues for possible automatic
bandwidth choice balancing the trade-off between variance and squared-bias
effects. However, the constants are unknown in practice which could seem to
be a serious drawback for practical purposes. This general problem is
well-known in classical nonparametric statistics, but in our functional
context this question gets even more intricate because of the rather
complicated expression of $M_{0}$, $M_{1}$ and $M_{2}$. An appealing way to
attack the problem is to use bootstrap ideas. In Section \ref{FuncWildBoot}
we propose a track for building a functional version of the so-called wild
bootstrap. We will show in Section \ref{FuncWildBootSim}, through some
simulated examples, how this functional wild bootstrap procedure works on
finite sample sizes for choosing automatically an optimal bandwidth. A case
study, based on spectrometric functional data coming from the food industry,
will be shortly presented in Section \ref{FuncWildBootSpect}.

At this stage it is worth noting that we have no asymptotic support for this
functional bootstrapping procedure. This open question will be one of the
main point discussed in the concluding Section \ref{Conc}.


\subsection{A Functional version of the wild bootstrap}

\label{FuncWildBoot} 
Basically, when using bootstrapping techniques one expects to approximate
directly the distribution of the error of estimation without having to
estimate the leading terms involved in some asymptotic expansion of this
error. In standard finite-dimensional problems (that is, when the variable $%
\mathcal{X}$ is valued in $\mathbb{R}^{p}$), a so-called wild bootstrap has
been constructed for approximating the distribution of the error of
estimation in kernel nonparametric regression. We refer to H\"{a}rdle (1989)
and H\"{a}rdle \& Marron (1991) for a previous presentation of the wild
bootstrap in nonparameric regression. A selected set of additional
references would include Mammen (2000) for the description of the state of
art on nonparametric bootstrapping, Mammen (1993) for a large study of wild
bootstrap, and H\"{a}rdle, Huet \& Jolivet (1995) for specific advances on
wild bootstrap in (finite-dimensional) nonparametric regression setting.

The main interest of this kind of bootstrap relies on a resampling procedure
of the residuals which makes it easily adaptable to our functional setting.
Precisely, an adaptation to our functional setting could be the following
functional wild bootstrap procedure:

\begin{itemize}
\item[i)] Given the estimate $\hat{r}$ constructed with a bandwidth $h$,
compute the residuals $\hat{\epsilon}_{i}=y_{i}-\hat{r}(\mathcal{X}_{i}),$
and construct a sequence of bootstrapped residuals such that each $\epsilon
_{i}^{\ast}$ is drawn from a distribution $G_{i}^{\ast}$ which is the sum of
two Dirac distributions : 
\begin{equation*}
G_{i}^{\ast}=\frac{5+\sqrt{5}}{10}\delta_{\frac{\hat{\epsilon}_{i}(1-\sqrt {5%
})}{2}}+(\frac{5-\sqrt{5}}{10})\delta_{\frac{\hat{\epsilon}_{i}(1+\sqrt{5})}{%
2}}.
\end{equation*}
Such a distribution ensures that the first three moments of the bootstrapped
residuals are respectively $0$, $\hat{\epsilon}_{i}^{2}$ and $\hat{\epsilon }%
_{i}^{3}$ (see H\"{a}rdle \& Marron, 1991, for details).

\item[ii)] Given the bootstrapped residuals $\epsilon_{i}^{\ast}$, and using
a new kernel estimate $\tilde{r}$ which is defined as $\hat{r}$ but by using
another bandwidth $g$, construct a bootstrapped sample $(\mathcal{X}%
_{i}^{\ast},Y_{i}^{\ast})$ by putting 
\begin{equation*}
\mathcal{X}_{i}^{\ast}=\mathcal{X}_{i}\ {\mbox{ and }}\ Y_{i}^{\ast}=\tilde {%
r}(\mathcal{X}_{i})+\epsilon_{i}^{\ast}.
\end{equation*}

\item[iii)] Given the bootstrapped sample $(\mathcal{X}_{i}^{\ast},Y_{i}^{%
\ast})$, compute the kernel estimate $\hat{r}^{\ast}(\chi)$ which is defined
as $\hat{r}$ (with the same bandwidth $h$) but using the bootstrapped sample 
$(\mathcal{X}_{i}^{\ast},Y_{i}^{\ast})$ instead of the previous sample $(%
\mathcal{X}_{i},Y_{i})$.\newline
\end{itemize}

\noindent We suggest to repeat several times (let say $N_{B}$ times) this
bootstrap procedure, and to use the empirical distribution of $\hat{r}%
^{*}(\chi) - \tilde{r}(\chi)$ for bandwidth selection purpose. Precisely,
the bootstrapped bandwidth is defined as follows:

\begin{definition}
\label{DefBB} Given $N_{B}$ replications of the above described
bootstrapping scheme, and given a fixed set H of bandwidths, the
bootstrapped bandwidth $h^{\ast}$ is defined by: 
\begin{equation*}
h^{\ast}=h^{\ast}(\chi)=\arg\min_{h\in H}\left( \frac{1}{N_{B}}\sum
_{b=1}^{N_{B}}(\hat{r}^{\ast}(\chi)-\tilde{r}(\chi))^{2}\right) .
\end{equation*}
\end{definition}

\noindent Of course, this procedure has still to be validated theoretically
(see discussion in Section \ref{Conc}), but we will see in the next Sections %
\ref{FuncWildBootSim} and \ref{FuncWildBootSpect} how it behaves both on
simulated on and real data samples.



\subsection{Some simulations}

\label{FuncWildBootSim} 

The aim of this section is to look at how the automatic bootstrapped
bandwidth constructed in Definition \ref{DefBB} behaves on simulated
samples. We construct random curves in the following way: 
\begin{equation*}
\mathcal{X}(t)= sin(\omega t) + (a+2\pi)t + b, \ t \in(-1,+1),
\end{equation*}
\noindent where $a$ and $b$ (respectively $\omega$) are r.r.v. drawn from a
uniform distribution on $(0,1)$ (respectively on $(0,2\pi)$). Some of these
curves are presented in Figure \ref{curves} below.


\begin{center}
\begin{figure}[h]
\centerline{\epsfig{figure=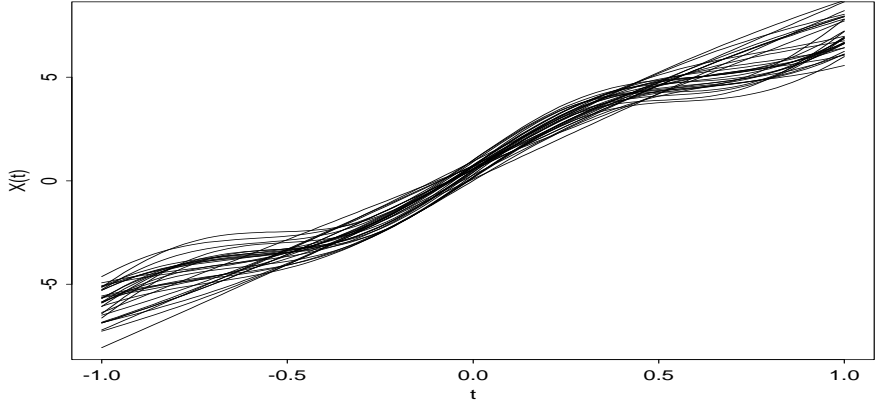,width=12cm, height=6cm}}
\caption{A sample of 30 simulated curves}
\label{curves}
\end{figure}
\end{center}

\noindent The real response is simulated according to the following
regression relation: 
\begin{equation*}
Y= r(\mathcal{X}) + \epsilon,
\end{equation*}
\noindent where 
\begin{equation*}
r(\mathcal{X})=\int_{-1}^{+1} |\mathcal{X}^{\prime}(t)|(1-cos(\pi t)) dt,
\end{equation*}
\noindent and where $\epsilon$ is drawn from a ${\mathcal{N}}(0,2)$
distribution.

Our experience is based on the following lines. For each experience, we
simulated two samples: a sample of size $n=100$ on which all the estimates
are computed and a testing sample of size $n=50$ which is used to look at
the behaviour of our method. Also, for each experience, the number of
bootstrap replications was taken to be $N_{B}=100$. Other values for $J$ and 
$N_{B}$ were also tried without changing the main conclusions. To improve
the speed of our algorithm, the bandwidth $h$ is assumed to belong to some
grid in terms of nearest neighbours, that is 
\begin{equation}  \label{Hsim}
h=h(\chi) \in\{h_{1}, \ldots, h_{32}\}=H,
\end{equation}
\noindent where $h_{k}$ is the radius of the ball of center $\chi$ and
containing exactly $k$ among the curves data $\mathcal{X}_{1}, \ldots 
\mathcal{X}_{100}$. Concerning the other parameters of our study, the kernel
function $K$ was chosen to be $K(u)=1-u^{2}, \ u \in(0,1)$ and the norm $%
||.||$ was taken to be the $L_{2}$ one between the first order derivatives
of the curves.

We computed, for the $32$ different values of $h$, the average (over the $%
\chi$'s belonging to the second testing data sample) of the true error $(%
\hat{r}(\chi)-r(\chi))^{2}$ and of its bootstrap approximation $(\hat
{r}%
^{\ast}(\chi)-\tilde{r}(\chi))^{2}$. Finally, this Monte Carlo scheme was
repeated $J=100$ times and the results are reported in Figure \ref%
{Compar-Errors} (only $25$ among the $100$ curves are presented to make the
plot clearer).

\begin{center}
\begin{figure}[h]
\centerline{\epsfig{figure=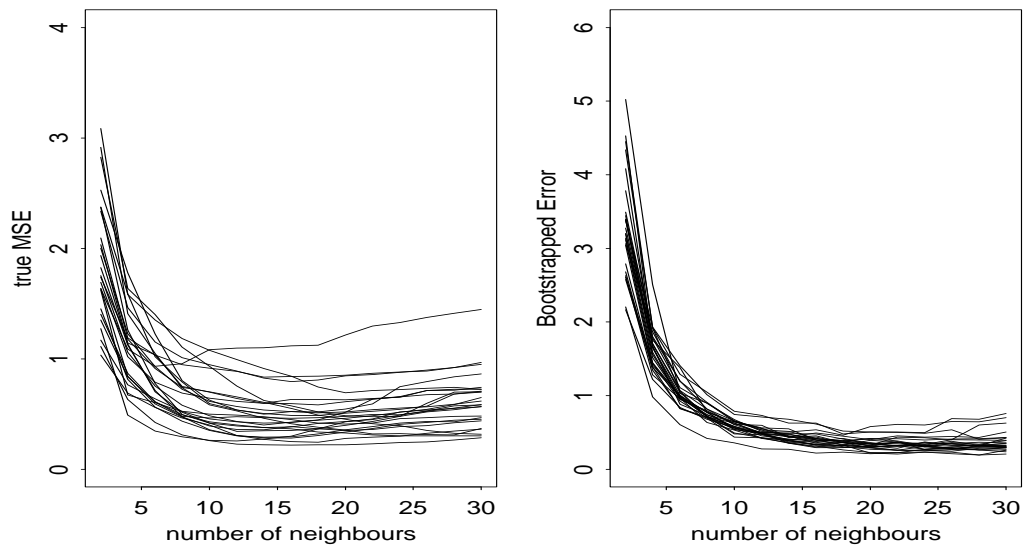,width=14cm, height=8cm}}
\caption{Simulations: True (left) and Bootstrap (right) Errors as functions
of h}
\label{Compar-Errors}
\end{figure}
\end{center}

\begin{center}
\begin{figure}[h]
\centerline{\epsfig{figure=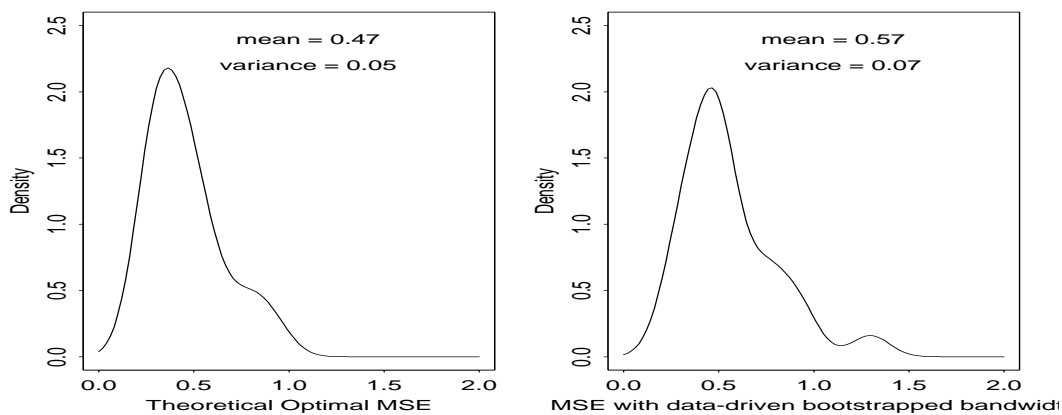,width=14cm, height=6cm}}
\caption{Simulations: MSE with optimal (left) and bootstrapped (right)
bandwidth}
\label{BandwChoiceResults}
\end{figure}
\end{center}

It appears clearly that both the theoretical quadratic loss and its
data-driven bootstrapped version have the same convex shape. This convex
shape is directly linked with the asymptotic expansion given in Theorem \ref%
{TheoremMSE} before: large values of $h$ give high bias, while small values
of $h$ lead to high variance. These results are quite promising in the sense
that the similarity of the shapes of both sets of curves presented in Figure %
\ref{Compar-Errors} let us expect that the bootstrapped bandwidths will be
closed from the optimal ones. To check that point , we computed the
theroretical minimal quadratic loss (that is, the error obtained by using
the best bandwidth) and we compared it with the error obtained by using the
boostrapped bandwidth $h^{\ast}$. This was done for each among the $J=100$
experiences, and the results are reported in Figure \ref{BandwChoiceResults}
which gives mean, variance and density estimates of these two errors.
Undoubtedly, these results show the good behaviour (at least on this
example) of the bootstrapping method as an automatic bandwidth selection
procedure. Of course, as discussed in Section \ref{Conc}, theoretical
support for this functional bootstrap bandwidth selection rule is still an
open question.

\subsection{A real data chemometric application}

\label{FuncWildBootSpect} 

\begin{center}
\begin{figure}[h]
\centerline{\epsfig{figure=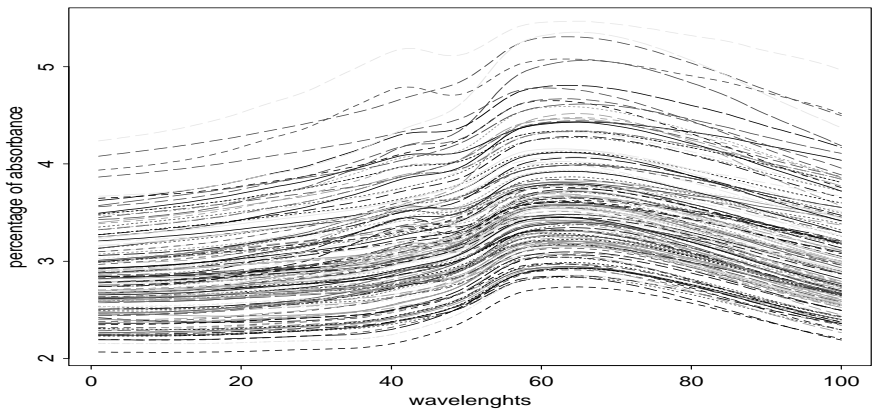,width=12cm, height=6cm}}
\caption{Spectrometric curves data}
\label{SpectroCurves}
\end{figure}
\end{center}

Let us now quickly show how our procedure is working on real data. These
data contain of $215$ spectra of light absorbance $(Z_{i},\ i=1,\ldots215)$
as functions of the wavelength, and observed on finely chopped pieces of
meat. We present in Figure \ref{SpectroCurves} the plots of the $215$
spectra. \newline

For each spectral cruve corresponds some real response $Y_{i}$ which is the
percentage of fatness, and our aim is to study the regression relation
existing between the real variable $Y$ and the functional one $Z$. These
data have been widely studied and, inspired by previous studies (see Ferraty
\& Vieu, 2006) we decide to apply the functional kernel methodology to the
curves $\mathcal{X}=Z^{\prime\prime}$, and by taking as norm $||.||$ between
curves the usual $L_{2}$ norm between the second derivatives of the spectra.
The kernel function $K$ was chosen to be $K(u)=1-u^{2},\ u\in(0,1)$. Along
our study we splitted the data into two subsamples. A first subsample of
size $n=165$ from which our estimates are computed, and a testing sample of
size $50$ on which they are applied.

In a first attempt, we used the automatic bootstrapping bandwidth selection
rule, where $H$ was defined as in (\ref{Hsim}). We present in Figure \ref%
{ErrorBootSpecCurves} the shape of the Bootsrapped Mean Square Error as a
function of the number of neighbours (and thus, as function of the
bandwidth).

\begin{center}
\begin{figure}[h]
\centerline{\epsfig{figure=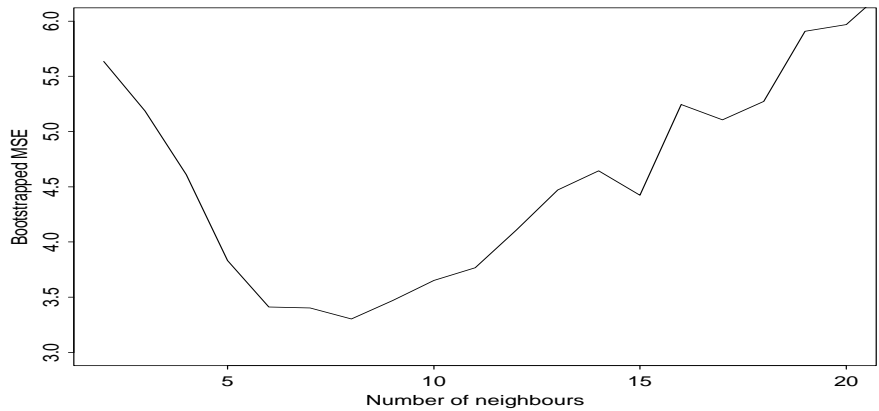,width=12cm, height=6cm}}
\caption{Spectrometric data: Bootstrapped Errors as function of the
bandwidth }
\label{ErrorBootSpecCurves}
\end{figure}
\end{center}

\noindent The same convex form as for the simulated data appears. This form
matches the theoretical results obtained in Section \ref{MSE}, with high
bias for large values of $h$ and high variance for small bandwidths. These
bootstrapped errors are, in this example, minimal for the value $k=8$. That
means that, for each new curve $\chi$ to be predicted, the data-driven
bootstrapped bandwidth $h^{\ast}(\chi)$ is such that there are exactly $8$
curves-data which are falling inside of the ball of radius $h^{\ast}(\chi)$.

\begin{center}
\begin{figure}[h]
\centerline{\epsfig{figure=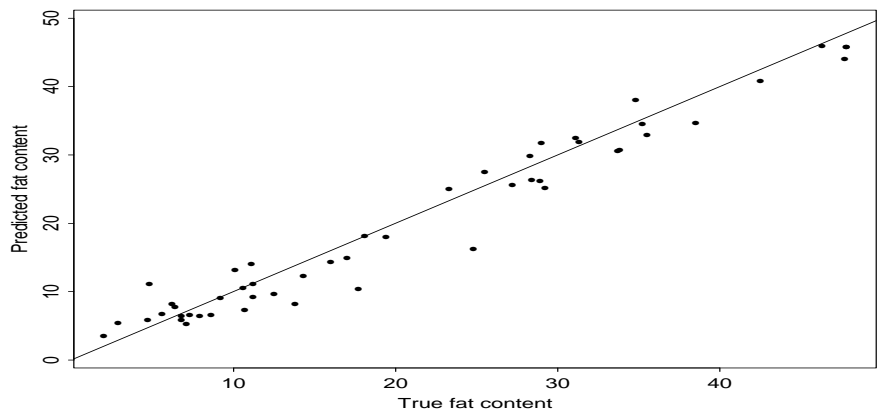,width=12cm, height=6cm}}
\caption{Spectrometric data: Predicted values on the testing sample}
\label{SpectroPredictions}
\end{figure}
\end{center}

These bandwidths lead to completely automatic data-driven fat contents
prediction. For instance, we present in Figure \ref{SpectroPredictions} the
fat content predictions for the $50$ spectra in our testing sample. In order
to highlight the nice behaviour of our prediction algorithm, Figure \ref%
{SpectroPredictions} plots the predicted values as functions of the true
ones.

\subsection{About implementation of the method}

\label{Splus}

The implementation of the method was performed by using the \emph{Splus}
routine \textit{funopare.kernel} which is included in the package \textit{%
npfda}. This package will go with the monograph by Ferraty \& Vieu (2006).
This \textit{Splus} package (as well as a similar \textit{R} package) and
the spectrometric dataset (as well as several other curves datasets) will be
put in free access on line in the next future. By that time, programs and
data are available on request. 

\section{Conclusions and open problems}

\label{Conc}


This paper completes the recent advances existing in kernel nonparametric
regression for functional data, by giving not only the rates of convergence
of the estimates but also the exact expressions of the constant terms
involved in these rates. These theoretical results deal with mean squared
errors evaluations and asymptotic normality results. As explained in Section %
\ref{ComputFeatures}, these new results open interesting perspectives for
applications, such as for instance data-driven automatic bandwidth selection
and confidence band construction.

We support the idea that bootstrap methods offers interesting perspectives
for the functional context. One of them is illustrated by our functional
version of the wild bootstrap for selecting the smoothing parameter. We have
observed nice results for bandwidth selection on some simulated and real
datasets. It should be pointed out that theoretical support for this
functional wild bootstrap bandwidth selection rule remains an open problem.
Our guess is that it should be possible to extend to functional variables
some results stated in finite dimension (for instance those in H\"{a}rdle \&
Bowman, 1987), but this has still to be proved.

Another direct application of our result concerns the construction of
confidence bands. Once again, the problem of estimating the constants
involved in the asymptotic normal distribution can be attacked by the wild
bootstrap track described before. One possible way for that would be to try
to extend standard finite-dimensional knowledge (see for instance H\"{a}rdle
\& Marron, 1991, or H\"{a}rdle, Huet \& Jolivet, 1995) to
infinite-dimensional variables.



\section{Appendix: Proofs}

\label{Appendix}

In the following, given some $\mathbb{R}$-valued random variable $U$, $P^{U}$
will stand for the probability measure induced by $U.$ To make its treatment
easier, the kernel estimate $\hat{r}$ will be decomposed as follows:

\begin{equation*}
\widehat{r}\left( \chi\right) =\frac{\widehat{g}\left( \chi\right) }{%
\widehat{f}\left( \chi\right) },
\end{equation*}
where 
\begin{equation*}
\widehat{g}\left( \chi\right) =\frac{1}{nF\left( h\right) }%
\sum_{k=1}^{n}Y_{k}K\left( \frac{\left\Vert \mathcal{X}_{k}-\chi\right\Vert 
}{h}\right)
\end{equation*}
and 
\begin{equation*}
\widehat{f}\left( \chi\right) =\frac{1}{nF\left( h\right) }%
\sum_{k=1}^{n}K\left( \frac{\left\Vert \mathcal{X}_{k}-\chi\right\Vert }{h}%
\right) .
\end{equation*}

\subsection{Proof of Theorem 1}

The proof is split into two parts: computations of the bias and of the
variance of the estimate. Each part is decomposed in technical lemmas that
will be proved in Section \ref{Appfin}.

\begin{itemize}
\item[-] \textbf{Bias term: proof of (\ref{B})}. Let us write the following
decomposition. 
\begin{equation}
E\widehat{r}\left( \chi\right) =\frac{E\widehat{g}\left( \chi\right) }{E%
\widehat{f}\left( \chi\right) }+\frac{A_{1}}{\left( E\widehat{f}\left(
\chi\right) \right) ^{2}}+\frac{A_{2}}{\left( E\widehat{f}\left( \chi\right)
\right) ^{2}},  \label{Bias}
\end{equation}
with%
\begin{equation}
A_{1}=E\left[ \widehat{g}\left( \chi\right) \left( \widehat{f}\left(
\chi\right) -E\widehat{f}\left( \chi\right) \right) \right]  \label{A1}
\end{equation}
and%
\begin{equation}
A_{2}=E\left[ \left( \widehat{f}\left( \chi\right) -E\widehat{f}\left(
\chi\right) \right) ^{2}\widehat{r}\left( \chi\right) \right] .  \label{A2}
\end{equation}

The first step of the proof consists in rewritting the first term in right
hand side of the decomposition (\ref{Bias}) in the following way:

\begin{lemma}
\label{lemma-bias-1} We have: \label{B2} 
\begin{equation*}
\frac{E\widehat{g}\left( \chi\right) }{E\widehat{f}\left( \chi\right) }%
-r\left( \chi\right) =h\varphi^{\prime}\left( 0\right) I + o\left( h\right) ,
\end{equation*}
where 
\begin{equation*}
I=\frac{\int_{0}^{1}tK\left( t\right) dP^{\left\Vert \mathcal{X}%
-\chi\right\Vert /h}\left( t\right) }{\int_{0}^{1}K\left( t\right)
dP^{\left\Vert \mathcal{X}-\chi\right\Vert /h}\left( t\right) }.
\end{equation*}
\end{lemma}

In a second attempt the next lemma will provide the constant term involved
in this bias expression and its limit.

\begin{lemma}
\label{B3}We have : 
\begin{equation*}
I=\frac{K\left( 1\right) -\int_{0}^{1}\left( sK\left( s\right) \right)
^{\prime}\tau_{h}\left( s\right) ds}{K\left( 1\right)
-\int_{0}^{1}K^{\prime}\left( s\right) \tau_{h}\left( s\right) ds}\underset{%
n\rightarrow\infty}{\longrightarrow} \frac{M_{0}}{M_{1}}.
\end{equation*}
\end{lemma}

Finally, to finish this proof it suffices to prove that both last terms at
right hand side of (\ref{Bias}) are neglectible. This is done in next lemma.

\begin{lemma}
\label{B7} We have: 
\begin{equation*}
A_{1}\ =\ O\left( \left( nF\left( h\right) \right) ^{-1}\right) \mbox{ and }
A_{2}\ =\ O\left( \left( nF\left( h\right) \right) ^{-1}\right) .
\end{equation*}
\end{lemma}

\noindent So the proof of (\ref{B}) is complete.


\item[-] \textbf{Variance term: proof of (\ref{V})}. The starting point of
the proof is the following decomposition. This decomposition has been
obtained in earlier work by Collomb (1976) (see also Sarda \& Vieu (2000))
in the finite-dimensional case, but since the proof is only using analytic
arguments about Taylor expansion of the function $1/z$ around $0$, it
extends obviously to our functional setting:

\begin{align}
Var\left( \widehat{r}\left( \chi\right) \right) & =\frac{Var\widehat {g}%
\left( \chi\right) }{\left( E\widehat{f}\left( \chi\right) \right) ^{2}}-4%
\frac{E\widehat{g}\left( \chi\right) Cov\left( \widehat{g}\left( \chi\right)
,\widehat{f}\left( \chi\right) \right) }{\left( E\widehat {f}\left(
\chi\right) \right) ^{3}}  \label{Decvar} \\
& +3Var\widehat{f}\left( \chi\right) \frac{\left( E\widehat{g}\left(
\chi\right) \right) ^{2}}{\left( E\widehat{f}\left( \chi\right) \right) ^{4}}%
+o\left( \frac{1}{nF\left( h\right) }\right) .  \notag
\end{align}

Finally, the result (\ref{V}) will follow directly from this decomposition
together with both following lemmas.

\begin{lemma}
\label{B4bis} We have successively: 
\begin{align*}
E\widehat{f}\left( \chi\right) & \rightarrow K\left( 1\right) -\int
_{0}^{1}K^{\prime}\left( s\right) \tau_{0}\left( s\right) ds=M_{1}, \\
E\widehat{g}\left( \chi\right) & \rightarrow r\left( \chi\right) \left(
K\left( 1\right) -\int_{0}^{1}K^{\prime}\left( s\right) \tau_{0}\left(
s\right) ds\right) =r\left( \chi\right) M_{1}.
\end{align*}
\end{lemma}

\begin{lemma}
\label{B5} We have successively: 
\begin{align*}
\left( Var\widehat{f}\left( \chi\right) \right) & =\frac{M_{2}}{nF\left(
h\right) }\left( 1+o\left( 1\right) \right) , \\
\left( Var\widehat{g}\left( \chi\right) \right) & =\left(
\sigma_{\varepsilon}^{2} +r^{2}\left( \chi\right) \right) \frac{M_{2}}{%
nF\left( h\right) }\left( 1+o\left( 1\right) \right) , \\
Cov\left( \widehat{g}\left( \chi\right) ,\widehat{f}\left( \chi\right)
\right) & =r\left( \chi\right) \frac{M_{2}}{nF\left( h\right) }\left(
1+o\left( 1\right) \right) .
\end{align*}
\end{lemma}
\end{itemize}

\subsection{Proof of Theorem 2}

The following lemma states a preliminary pointwise limiting distribution
result. This lemma stems from the bias and variance expressions obtained
along Theorem 1; it will be proved in the next subsection.

\begin{lemma}
\label{TCL1} We have: 
\begin{equation}  \label{equ-tcl}
\sqrt{nF\left( h\right) }\left( \widehat{r}\left( \chi\right) -r\left(
\chi\right) -B_{n}\right) \frac{M_{1}}{\sqrt {\sigma_{\varepsilon}^{2} M_{2}}
}\hookrightarrow\mathcal{N}\left( 0,1\right) .
\end{equation}
\end{lemma}

Because of standard Glivenko-Cantelli type results, we have 
\begin{equation*}
\frac{\widehat{F}\left( h\right) }{F\left( h\right) }\overset{\mathbb{P}}{%
\rightarrow}1,
\end{equation*}
and this is enough, combined with the result of Lemma \ref{TCL1}, to get the
conclusion of Theorem 2.

\subsection{Proofs of technical lemmas}

\label{Appfin}

\begin{itemize}
\item[-] \textbf{Proof of Lemma \ref{lemma-bias-1}:} To calculate 
\begin{equation}
\frac{E\widehat{g}\left( \chi\right) }{E\widehat{f}\left( \chi\right) }%
-r\left( \chi\right) =\frac{E\left[ \left( Y-r\left( \chi\right) \right)
K\left( \frac{\left\Vert \mathcal{X}-\chi\right\Vert }{h}\right) \right] }{%
EK\left( \frac{\left\Vert \mathcal{X}-\chi\right\Vert }{h}\right) },
\label{Pro1}
\end{equation}
note first that 
\begin{align*}
E\left[ \left( Y-r\left( \chi\right) \right) K\left( \frac{\left\Vert 
\mathcal{X}-\chi\right\Vert }{h}\right) \right] & =E\left[ \left( r\left( 
\mathcal{X}\right) -r\left( \chi\right) \right) K\left( \frac{\left\Vert 
\mathcal{X}-\chi\right\Vert }{h}\right) \right] \\
& =E\left[ \varphi\left( \left\Vert \mathcal{X}-\chi\right\Vert \right)
K\left( \frac{\left\Vert \mathcal{X}-\chi\right\Vert }{h}\right) \right] .
\end{align*}
Moreover, it comes: 
\begin{align*}
& E\left[ \varphi\left( \left\Vert \mathcal{X}-\chi\right\Vert \right)
K\left( \frac{\left\Vert \mathcal{X}-\chi\right\Vert }{h}\right) \right] \\
& =\int\varphi\left( t\right) K\left( \frac{t}{h}\right) dP^{\left\Vert 
\mathcal{X}-\chi\right\Vert }\left( t\right) \\
& =\int\varphi\left( ht\right) K\left( t\right) dP^{\left\Vert \mathcal{X}%
-\chi\right\Vert /h}\left( t\right) \\
& =h\varphi^{\prime}\left( 0\right) \int tK\left( t\right) dP^{\left\Vert 
\mathcal{X}-\chi\right\Vert /h}\left( t\right) +o\left( h\right) ,
\end{align*}

the last line coming from the first order Taylor's expansion for $\varphi$
around $0$. For the denumerator in (\ref{Pro1}) we have%
\begin{equation*}
E\left[ K\left( \left\Vert \mathcal{X}-\chi\right\Vert /h\right) \right]
=\int K\left( t\right) dP^{\left\Vert \mathcal{X}-\chi\right\Vert /h}\left(
t\right) .
\end{equation*}

\noindent Finally, it appears clearly that the first order bias term is $%
h\varphi^{\prime}\left( 0\right) I$.

\item[-] \textbf{Proof of Lemma \ref{B3}:} We note that%
\begin{equation*}
tK\left( t\right) =K\left( 1\right) -\int_{t}^{1}\left( sK\left( s\right)
\right) ^{\prime}ds.
\end{equation*}
Applying Fubini's Theorem we get%
\begin{align*}
\int_{0}^{1}tK\left( t\right) dP^{\left\Vert \mathcal{X}-\chi\right\Vert
/h}\left( t\right) & =K\left( 1\right) F\left( h\right) -\int_{0}^{1}\left(
\int_{t}^{1}\left( sK\left( s\right) \right) ^{\prime }ds\right)
dP^{\left\Vert \mathcal{X}-\chi\right\Vert /h}\left( t\right) \\
& =K\left( 1\right) F\left( h\right) -\int_{0}^{1}\left( sK\left( s\right)
\right) ^{\prime}F\left( hs\right) ds.
\end{align*}
Similarly, we have 
\begin{equation}  \label{Pro2}
\int_{0}^{1}K\left( t\right) dP^{\left\Vert \mathcal{X}-\chi\right\Vert
/h}\left( t\right) =K\left( 1\right) F\left( h\right)
-\int_{0}^{1}K^{\prime}\left( s\right) F\left( hs\right) ds.
\end{equation}
So the proof of this lemma is finished by applying the Lebesgue's dominated
convergence theorem since the denumerator may be easily bounded above by $%
K\left( 1\right) >0$ ($K$ being decreasing).

\item[-] \textbf{Proof of Lemma \ref{B4bis}}: The first assertion follows
directly from (\ref{Pro2}), while the second one can be proved similarly
according to the following lines: 
\begin{align*}
EYK\left( \frac{\left\Vert \mathcal{X}-\chi\right\Vert }{h}\right) & =E\left[
E\left( Y|\mathcal{X}\right) K\left( \frac{\left\Vert \mathcal{X}%
-\chi\right\Vert }{h}\right) \right] \\
& =\left( r\left( \chi\right) +o(1)\right) E\left[ K\left( \frac{\left\Vert 
\mathcal{X}-\chi\right\Vert }{h}\right) \right] .
\end{align*}

\item[-] \textbf{Proof of Lemma \ref{B5}}: We write the variance of $%
\widehat{f}\left( \chi\right) $ as:

\begin{equation*}
\left( Var\widehat{f}\left( \chi\right) \right) =\frac{1}{nF^{2}\left(
h\right) }\left[ EK^{2}\left( \frac{\left\Vert \mathcal{X}-\chi\right\Vert }{%
h}\right) -\left( EK\left( \frac{\left\Vert \mathcal{X}-\chi\right\Vert }{h}%
\right) \right) ^{2}\right] ,
\end{equation*}
and note that, as for getting (\ref{Pro2}), it holds: 
\begin{equation}  \label{EK2}
EK^{2}\left( \frac{\left\Vert \mathcal{X}-\chi\right\Vert }{h}\right)
=F\left( h\right) \left( K^{2}\left( 1\right) -\int_{0}^{1}\left(
K^{2}\right) ^{\prime}\left( s\right) \tau_{h}\left( s\right) ds\right) .
\end{equation}
The first assertion of lemma \ref{B4bis} gives 
\begin{equation*}
\left( EK\left( \frac{\left\Vert \mathcal{X}-\chi\right\Vert }{h}\right)
\right) ^{2}=O\left( F^{2}\left( h\right) \right) .
\end{equation*}
At last 
\begin{equation}
\left( Var\widehat{f}\left( \chi\right) \right) \sim(n\,F\left( h\right)
)^{-1}\left( K^{2}\left( 1\right) -\int_{0}^{1}\left( K^{2}\right)
^{\prime}\left( s\right) \tau_{0}\left( s\right) ds\right) ,  \label{Pro3}
\end{equation}
which finishes the proof of the first assertion of our lemma.

The same steps can be followed to prove the second assertion. We write 
\begin{equation*}
\left( Var\widehat{g}\left( \chi\right) \right) =\frac{1}{nF^{2}\left(
h\right) }\left[ EY^{2}K^{2}\left( \frac{\left\Vert \mathcal{X}%
-\chi\right\Vert }{h}\right) -\left( EYK\left( \frac{\left\Vert \mathcal{X}%
-\chi\right\Vert }{h}\right) \right) ^{2}\right] .
\end{equation*}
The second term at right hand side of this expression is treated directly by
using the second assertion of Lemma \ref{B4bis}, while the first one is
treated as follows by conditioning on $\mathcal{X}$: 
\begin{equation*}
EY^{2}K^{2}\left( \frac{\left\Vert \mathcal{X}-\chi\right\Vert }{h}\right) \
=\ Er^{2}\left( \mathcal{X}\right) K^{2}\left( \frac{\left\Vert \mathcal{X}%
-\chi\right\Vert }{h}\right) + E\sigma_{\varepsilon}^{2}(\mathcal{X}%
)K^{2}\left( \frac{\left\Vert \mathcal{X}-\chi\right\Vert }{h}\right) .
\end{equation*}
The continuity of $r^{2}$ and of $\sigma^{2}_{\varepsilon}(.)$ insure that 
\begin{equation*}
Var\widehat{g}\left( \chi\right) =\frac{1}{nF^{2}\left( h\right) }\left(
\sigma_{\varepsilon}^{2}+r^{2}\left( \chi\right) \right) EK^{2}\left( \frac{%
\left\Vert \mathcal{X}-\chi\right\Vert }{h}\right) \left( 1+o\left( 1\right)
\right) .
\end{equation*}
\noindent Combining this result with (\ref{EK2}) allows to finish the proof
of the second assertion of our lemma. Let us deal now with the covariance
term :%
\begin{align}
& Cov\left( \widehat{g}\left( \chi\right) ,\widehat{f}\left( \chi\right)
\right)  \notag \\
& =\frac{1}{nF^{2}\left( h\right) }\left[ EYK^{2}\left( \frac{\left\Vert 
\mathcal{X}-\chi\right\Vert }{h}\right) -EK\left( \frac{\left\Vert \mathcal{X%
}-\chi\right\Vert }{h}\right) EYK\left( \frac{\left\Vert \mathcal{X}%
-\chi\right\Vert }{h}\right) \right] .  \notag
\end{align}
The last two terms were computed before, while the first one is treated by
conditioning on $\mathcal{X}$ and using continuity of $r$: 
\begin{equation*}
EYK^{2}\left( \frac{\left\Vert \mathcal{X}-\chi\right\Vert }{h}\right)
=\left( r\left( \chi\right) +o(1) \right) EK^{2}\left( \frac{\left\Vert 
\mathcal{X}-\chi\right\Vert }{h}\right) .
\end{equation*}
The proof of this lemma is now finished.

\item[-] \textbf{Proof of Lemma \ref{B7}}: Both assertions of this lemma are
direct consequences of Lemmas \ref{B4bis} and \ref{B5}.

\item[-] \textbf{Proof of Lemma \ref{TCL1}:} On one hand, (\ref{Bias}) and
Lemma \ref{B7} allows us to get 
\begin{equation*}
\widehat{r}(\chi)-E\widehat{r}(\chi)\ =\ \frac{\widehat{g}(\chi)}{\widehat {f%
}(\chi)}-\frac{E\widehat{g}(\chi)}{E\widehat{f}(\chi)}+o\left( \frac {1}{%
\sqrt{n\,F(h)}}\right) .
\end{equation*}
On the other hand, the following decomposition holds: 
\begin{equation*}
\frac{\widehat{g}(\chi)}{\widehat{f}(\chi)}-\frac{E\widehat{g}(\chi )}{E%
\widehat{f}(\chi)}=\frac{\left( \widehat{g}(\chi)-E\widehat{g}(\chi)\right) E%
\widehat{f}(\chi)+\left( E\widehat{f}(\chi)-\widehat{f}(\chi)\right) E%
\widehat{g}(\chi)}{\widehat{f}(\chi)\,E\widehat{f}(\chi)}.
\end{equation*}
Using Slutsky's theorem and Theorem \ref{TheoremMSE}, we get 
\begin{equation*}
\sqrt{\frac{n\,F(h)\,M_{1}^{2}}{\sigma_{\varepsilon}^{2}M_{2}}}\left( 
\widehat{r}(\chi)-E\widehat{r}(\chi)\right) \hookrightarrow\mathcal{N}\left(
0,1\right) ,
\end{equation*}
noting that $\widehat{r}(\chi)-E\widehat{r}(\chi)$ can be expressed as an
array of independent centered random variables (and the Central Limit
Theorem applies). Let us remark that 
\begin{equation*}
\widehat{r}(\chi)-E\widehat{r}(\chi)\ =\ \widehat{r}\left( \chi\right)
-r\left( \chi\right) -B_{n},
\end{equation*}
which achieves the proof of this lemma.
\end{itemize}

\subsection{Proof of Propositions \protect\ref{prop1} and \protect\ref{Pr2}}

\begin{itemize}
\item[-] \textbf{Proof of Proposition \ref{prop1}-$i$ and $ii$}: Obvious.

\item[-] \textbf{Proof of Proposition \ref{prop1}-$iii$}: We have: 
\begin{equation*}
\tau_{h}\left( s\right) =\exp\left( -\frac{C_{2}}{h^{p}}\left( \frac {1}{%
s^{p}}-1\right) \right) s\gamma.
\end{equation*}
If $s=1$ we have $\tau_{h}\left( 1\right) =1, \forall h$, while if $s=0$ we
have $\tau_{h}\left( 1\right) =\tau_{0}\left( 1\right) = F(0)=0, \forall h$.
To complete this proof it suffices to note that, for $s\in\left( 0,1\right) $%
, we have: $1/s^{p}-1>0,$ and so we have $\tau_{h}\left( s\right)
\rightarrow0$ as $h\rightarrow0$.

\item[-] \textbf{Proof of Propostion \ref{prop1}-$iv$}: For any $s>0$ and
any $h<1$, we have: 
\begin{equation*}
\tau_{h}\left( s\right) =\frac{\left\vert \ln h\right\vert }{\left\vert \ln
h+\ln s\right\vert }=\frac{\ln h}{\ln h+\ln s}=\frac{1}{1+\frac{\ln s}{\ln h}%
}
\end{equation*}
and so we have $\tau_{0}(s)=1, \forall s>0$. To complete this proof it
suffices to note that $\tau_{h}\left( 1\right) =\tau_{0}\left( 1\right) =
F(0)=0, \forall h$.

\item[-] \textbf{Proof of Propostion \ref{Pr2}-$i$}: By simple integration
by parts we arrive at 
\begin{equation*}
M_{0}=\int_{0}^{1}sK(s)T_{0}(s)ds,
\end{equation*}
\noindent where $T_{0}^{\prime}=\tau_{0}$. Because of $\mathbf{H2}$ and
because $\tau_{0}$ is non decreasing, there exists some nonempty interval $%
[a,b]\subset(0,1)$ such that both $K$ and $T_{0}$ do not vanish on $[a,b]$,
and therefore we arrive at: 
\begin{equation*}
M_{0}\geq\int_{a}^{b}sK(s)T_{0}(s)ds>0.
\end{equation*}

\item[-] \textbf{Proof of Propostion \ref{Pr2}-$i$}: Obvious.
\end{itemize}

\noindent{\small {\textbf{Acknowledgements.} The authors are grateful to all
the participants of the working group STAPH
(http://www.lsp.ups-tlse.fr/staph) on Functional Statistics at the
University Paul Sabatier of Toulouse for their stimulating and helpful
comments, with special thanks to Laurent Delsol for his cosnructie reading
of a earlier version of this paper. This work has been significantly
improved along the refereing procedure, and we would like to thank two
anonymous referees, an Associate Editor and the Editor for their comments
and suggestions.} }

{\small 
}

\end{document}